\documentclass{article}
\title{Existence of $C^{1,\alpha}$ Singular Solutions to Euler-Nernst-Planck-Poisson System on $\mathbb{R}^3$ with Free-Moving Charges}
\date{}
\author{Yiya Qiu\footnote{CAS Wu Wen-Tsun Key Labotatory of Mathematics and School of Mathematical Sciences, University of Science and Technology of China, Hefei, Anhui, 230026, PR China}, Lifeng Zhao\footnote{School of Mathematical Sciences, University of Science and Technology of China, Hefei, Anhui, 230026, PR China}}

\makeatletter
\@addtoreset{equation}{section}
\makeatother

\usepackage{amsmath}
\usepackage{amssymb}
\usepackage{amsthm}
\usepackage{mathtools}
\usepackage[a4paper,left=3cm,right=3cm,top=3cm,bottom=3cm]{geometry}
\usepackage{subeqnarray}

\newtheorem{thm}{Theorem}[section]
\newtheorem{lem}[thm]{Lemma}
\newtheorem{prop}[thm]{Proposition}

\theoremstyle{remark}
\newtheorem{rmk}{Remark}
\newcommand\al{\alpha} 
\newcommand\be{\beta}  
\newcommand\ga{\gamma}  
\newcommand\ve{\varepsilon}  
\newcommand\la{\lambda}    
\newcommand\w{\omega}
\newcommand\W{\Omega}
\newcommand\se{\theta}

\newcommand\vs{\vartheta}

\newcommand\ma{\mathcal}

\newcommand\les{\lesssim} 
\newcommand\f{\frac}  
\newcommand\p{\partial}  

\newcommand\td{\tilde} 
\newcommand\mc{\mathcal}

\begin{document}
\maketitle
\begin{abstract}
We construct a special $C^{1,\al}$ blow up solution to the three dimensional system modeling electro-hydrodynamics, which is a strongly coupled system of incompressible Euler equation and Nernst-Planck-Poisson equation. Our construction lies on the framework established in \cite{Elgindi2019a} and relies on a special solution to variant spherical Laplacian.
\end{abstract}
\bibliographystyle{plain}

\section{Introduction}
In this article, we consider Euler-Nernst-Planck-Poisson system with free-moving charges on $\mathbb{R}^3$
\begin{equation} \label{1.1} 
\begin{cases}
\p_t u+(u\cdot\nabla)u+\nabla p=\Delta\xi\nabla\xi,\\
\nabla\cdot u=0,\\
\p_tv+(u\cdot\nabla)v=0,\\
\Delta\xi=v,\\
(u,v,\xi)|_{t=0}=(u_0,v_0,\xi_0),
\end{cases}
\end{equation}
which is a special case of Navier-Stokes-Poisson-Nernst-Planck system
\begin{equation}\label{origin}
\begin{cases}
\p_tu+(u\cdot\nabla)u-\vs\Delta u+\nabla p=\bar{\ve}\Delta \xi \nabla \xi,\\
\p_t w_1+(u\cdot\nabla)w_1=\nabla\cdot(D_1\nabla w_1-\iota_1 w_1\nabla \xi),\\
\p_t w_2+(u\cdot\nabla)w_2=\nabla\cdot(D_2\nabla w_2+\iota_2 w_2\nabla \xi),\\
\nabla\cdot u=0,\\
\bar{\ve}\Delta\xi=w_1-w_2,\\
(u,w_1,w_2,\xi)|_{t=0}=(u_0,w_{1,0},w_{2,0},\xi_0).
\end{cases}
\end{equation}
In the system (\ref{origin}), $u:\mathbb{R}^3\times[0,\infty)\to\mathbb{R}^3$ is the velocity field; $p:\mathbb{R}^3\times[0,\infty)\to\mathbb{R}$ is the scalar pressure; $\xi:\mathbb{R}^3\times[0,\infty)\to\mathbb{R}$ is the electrostatic potential; $w_1$ and $w_2$ are the densities of binary diffuse negative and positive charges, respectively. $\vs$ is the fluid viscosity, $\bar{\ve}$ is the dielectric constant of the fluid. $D_1$, $D_2$, $\iota_1$, $\iota_2$ are the diffusion and mobility coefficients of the charges. In this article, we mainly consider the inviscid flow, namely, $\vs=0$, and we let $\bar{\ve}=1$ without loss of generality. Moreover, since diffusion in charges moving is a low order effect, which can be safely neglected, so we let $D_1=D_2=\iota_1=\iota_2=0$. At last we let $v=w_1-w_2$ and get the system (\ref{1.1}).

First introduced by Rubinstein in \cite{Rubinstein1990}, the system (\ref{origin}) is a well-established mathematical model for electro-hydrodynamics, coupling the Navier-Stokes of an incompressible fluid and the transported Poisson-Nernst-Planck equations.
It describes the dynamical behavior of incompressible fluid with binary diffuse positive and nagative charges (e.g., ions), which is relevant to understanding the behavior of different physical objects such as electrolyte solution, ion-exchangers, ion-selective membranes, and semiconductors. Over the years, (\ref{origin}) is well applied in many scientific areas such as electrical, fluid-mechanical and bio-chamical phenomena occuring in complex bio-hydred system, see \cite{Balbuena2004}, \cite{Enikov2000}, \cite{Enikov2005}, \cite{Shahinpoor2004} for various applications and \cite{Longaretti2008}, \cite{Jerome2008} for computational results.  For better understanding about the physical background and mathematical description, one would refer to \cite{Biler2000}, \cite{Bazant2004}, \cite{Probstein2003}, \cite{Chu2005} and  \cite{SCHMUCK2009}.

The mathematical analysis of (\ref{origin}) has been well studied by many experts for dozens of years, and their works can not be totally listed here. To our knowledge, the results date back to \cite{Jerome2002}, where the author established the local wellposedness of the system by semigroup theory. Then, in \cite{Jerome2009} and \cite{Fischer2016} the authors proved the existence of global weak solutions with block boundary.
The global existence for small data holds in \cite{Ryham2009} via linearization of a relative entropy functional.  Besides, global existence for mild solutions with small data in various function spaces such as Triebel–Lizorkin space and negative-order Besov space were obtained in \cite{Deng2011} and \cite{Zhao2010}. And global existence for large data in 2D was got by \cite{SCHMUCK2009} via energy laws, mass conservation, non-negativity and pointwise bounds. Recently in \cite{Bothe2014},  the authors used energy estimate and got exponential stability of steady states, and global existence for blocking and general selective boundary conditions was proved in \cite{Constantin2018} by exploiting Boltzmann states.

As for the blow up aspects, there has been several blowup criterions for (\ref{origin}) obtained in \cite{Zhao2016} and \cite{Liu2019}. However, the existence of blowup solutions still remains open, let alone the blowup dynamics for some appropriately chosen initial data. Indeed, it is rather challenging to construct blow up solutions even for Euler and Navier-Stokes. But recently, the progress concerning formation of self-similar singular solutions in fluid dynamics is flourishing. For example, in the area of incompressible fluid, Elgindi constructed the $C^{1,\al}$ singular solution to Euler equation in \cite{Elgindi2019}. At first, he established an approximation of Biot-Savart law naming fundamental model in a variant spherical coordinate. This model not only can be solved explicitly, but also wisely separated the solution into singular part and regular part. Then he performed linearization around the fundamental model and verified the coericivity of linearized operator in some weighted Sobolev space. At last with elliptic estimate and energy estimates, Elgindi finally established the formation of singular solutions to Euler in $\mathbb{R}^3$ with $C^{1,\al}$ data.
While in the field of compressible fluid, by means of Riemann variable and 3D Burgers profile, the finite-time blowup with smooth data of 3D compressible Euler was obtained in \cite{Buckmaster2019} without symmetry assumptions. Besides, in \cite{Merle2019}, the authors constructed global self-similar profile of 3D compressilble Euler by using classical phase portrait of nonliear ODEs. In addition, self-similar profile of 2D Burgers and the spectrum of linearzed operators were analyzed in \cite{Collot2018a} to get the singularity formation of 2D unsteady Prandtl's system.

In this article, following the framework of \cite{Elgindi2019} and \cite{Elgindi2019a} built by Elgindi, we prove the existence of finite-time blowup solutions to (\ref{1.1}).
\begin{thm} \label{thm}
There exists an $\al>0$, a divergence-free $u_0\in C^{1,\al}(\mathbb{R}^3)$ and a function $v_0\in C^\infty(\mathbb{R}^3)$ with initial vorticity $|\w_0(x)|\leq\f{C}{|x|^\al+1}$ and initial charges density $|v_0(x)|\leq\f{C}{|x|^{\f{5-\sqrt{21}}{2}}+1}$ for some constant $C>0$ so that the unique local solution $(u,v)\in (C_{x,t}^{1,\al}([0,1)\times\mathbb{R}^3)\times C_{x,t}^{\infty}([0,1)\times\mathbb{R}^3))$ to (\ref{1.1}) satisfies
$$\lim_{t\to1}\int_{0}^{t}\|\w\|_{L^\infty}ds=+\infty,\ \ \ \ \lim_{t\to1}\|v(t)\|_{L^\infty}=0.$$
\end{thm}
We shall explain our main ideas. First, the evolution of $v$ is a linear transport process with velocity of $u$,  so $v$ can not develop singularity by itself because $\|u\|_{L^\infty}$ doesn't blow up near the origin. 
So it is reasonable to construct $(u,v)$ where $v$ was kept small during the whole lifespan. In this way, the influence of $v$ on $u$ is under control once we make proper assumption on $v$.
So in our proof, the electrostatic effect is regarded as a small perturbation to incompressible fluid dynamics.
Next, $v$ has to be performed coordinate transformation following $u$, then we show the coericivity of linearized operator under the new coordinate.
After that, we analyze the contribution of $v$ to $u$, which is nonlocal because the effect is imposed by the electrostatic potential $\xi$ relating to $v$ by Poisson equation. To overcome the obstacle, with classical Strum-Liouville theory, we try to find a special explicit solution to the variant spherical Laplacian under new $(z,\se)$ coordinate (see Section 4 for details), where there is only one parameter under determined. At last, we follow \cite{Elgindi2019a} to derive the modulation law and perform the energy estimate, finishing the proof of Theorem 1.1.

\begin{rmk}
Here we only choose one special solution of spherical Laplacian, but there are still others can be used to construct the solutions for system, so the solution we find is not unique.
\end{rmk}

\begin{rmk}
This result is NOT a stability one. In our construction the perturbation of electrostatic effect doesn't occur in an open set. Indeed the profile of $v_0(x)$ is fixed, and what happens if we perturb the initial data is still unknown. To our sight, in order to do stability estimate, one must establish elliptic estimate to (\ref{ellp2}). However, we have to control $\|\sec\se\Pi\|_{L^2}$ appearing in the $L^2$ estimate, which seems hard to handle.
\end{rmk}

\begin{rmk}
In our construction, $u$ blows up in finite time, while $v$ goes to 0 near the  blow up time. A natural question is whether there exists a solution that $u$ and $v$ blow up at the same time. However, in this case $v$ can no longer be regareded as the perturbation of $u$ so that the contribution of $v$ to $u$ seems hard to control near the blowup time, because the nontrivial profile of $v$ may change the profile of $u$.
\end{rmk}

This aritcle is organized as follows:
In Section 2 we recall the notations and some lemmas used in \cite{Elgindi2019} and \cite{Elgindi2019a}. In Section 3 we follow the framework of \cite{Elgindi2019a} to perform the self-similar transformation and linearization. In Section 4 we discuss the coericivity of linearized operators. In Section 5 we state the elliptic estimate and construct the special solution of variant spherical Laplacian. In Section 6 we derive the modulation law. In Section 7 we perform the energy estimate so that main result follows.

\section{Preliminary}
We inherit the notations and some lemmas from \cite{Elgindi2019} and \cite{Elgindi2019a}.
\subsection{Notations}

$r$ will denote the two dimensional radial variable:
$$r=\sqrt{x_1^2+x_2^2}.$$
$\rho$ wil denote three dimensional variable:
$$\rho=\sqrt{r^2+x_3^2}.$$
And $\rho^\al$ is denoted by $R$:
$$R=\rho^\al.$$
We write $\zeta$:
$$\zeta=\arctan(\f{x_2}{x_1})$$ 
and $\se$:
$$\se=\arctan(\f{x_3}{r}).$$ 
$z$ will denote generally self-similar radial variable:
$$z=\f{R}{(1-(1+\mu)t)^{(1+\la)}},$$ where $\la$ and $\mu$ are small constants.
 Obviously, $\zeta,\se\in[0,\f{\pi}{2}]$, $z\in[0,+\infty)$. And $\se=0$ corresponds to the plane $x_3=0$ while $\se=\f{\pi}{2}$ corresponds to the $x_3$ axis. The main parameters we will use are $$0<\al\ll1,\ \ \ \ \ga=1+\f{\al}{10},\ \ \ \ \eta=\f{99}{100}.$$

Let $w$ be a weight function:
$$w=\f{(1+z)^2}{z^2},\ \ \ \ z\in[0,+\infty).$$

Besides, we define
$$\Gamma(\se)=(\sin\se\cos^2\se)^\f{\al}{3}$$ and 
\begin{equation}\label{K}
K(\se)=3\sin\se\cos^2\se.
\end{equation}

Define the differential operators:
$$D_\se(f)=\sin2\se\p_\se f,\ \ \\ \ D_R(f)=R\p_Rf.$$

Define the $L^2$ inner product and norm:$$(f,g)_{L^2_{\W}}=\int_\W fg,\ \ \ \ 
\|f\|_{L^2}=\sqrt{(f,f)_{L^2_\W}}.$$

Define $L^\infty$ norm:$$\|f\|_{L^\infty_\W}=\sup_{x\in\W}|f(x)|.$$

Define norm of space $C^{0,\be}$:
$$\|f\|_{C^\be_\W}=\sup_{x\in\W}|f|+\sup_{x\neq y}\f{|f(x)-f(y)|}{|x-y|^\be}.$$If $f\in C^1$ we say that $f\in C^{1,\be}$ when $\nabla f\in C^\be$.

Define weighted Sobolev norm $\mathcal{H}^k([0,\infty)\times[0,\pi/2])$:
\begin{equation}\label{Hk}
\|f\|^2_{\mathcal{H}^k}=\sum_{i=0}^{k}\|D^i_Rf\f{w}{\sin^{\eta/2}2\se}\|^2_{L^2}+
\sum_{0\leq i+j\leq k,i\geq 1}\|D^j_R D^i_\se f\f{w}{\sin^{\ga/2}2\se}\|^2_{L^2}.
\end{equation}

For any $f:[0,+\infty)\times[0,\f{\pi}{2}]\to \mathbb{R}$, we define non-local linear operator 
\begin{equation*}
L_{K}(f)(z)=\int_z^\infty\int_0^\f{\pi}{2}\f{f(z'){K}(\se)}{z'}dz'd\se,
\end{equation*}where $K$ is given in (\ref{K}).

Define the linearized operators $\mathcal{L}$, $\mc{L}_{F_*}$ and $\mc{L}^T_{F_*}$:
$$\mc{L}(f)=f+z\p_z f-2\f{f}{1+z},$$
$$\mc{L}_{F_*}(f)=f+z\p_z f-2\f{f}{1+z}-\f{2z\Gamma(\se)}{c(1+z)^2}L_{K}(f),$$
and$$\mc{L}_{F_*}^T(f)=\mc{L}_\Gamma(f)-\mathbb{P}(\f{3}{1+z}\sin2\se\p_\se f),$$
where$$\mathbb{P}(f)(z,\se)=f(z,\se)-\f{\Gamma(\se)}{c}\f{2z^2}{(1+z)^2}L_{K}(f)(0),\ \ \ c=\int_0^\f{\pi}{2}\Gamma(\se)K(\se)d\se.$$

\subsection{Some Useful Facts}
Let $$F_*(z,\se)=\f{\Gamma(\se)}{c}\f{4z}{(1+z)^2},$$ we have
\begin{equation*}
L_{K}(F_*)(z)=\f{4\al}{1+z}.
\end{equation*}

For any $f\in \mathcal{H}^k$, 
\begin{equation*}
L_{{K}}(z\p_z f)=\int_{z}^{\infty} \int_{0}^{\f{\pi}{2}} \p_p f(p,\se){K}(\se)d\se dp
=z\p_zL_{{K}}(f)(z),
\end{equation*}
and in particular,
\begin{equation} \label{zpz}
L_{{K}}(z\p_z f)(0)=-\int_{0}^{\f{\pi}{2}} f(0,\se){K}(\se)d\se.
\end{equation}

Besides, there is the following boundness for $L_K$:
\begin{prop}\emph{(\cite{Elgindi2019})}There exists a universal constant $C>0$ such that for all $f\in\ma{H}^k$ with $L_K(f)(0)=0$, we have
\begin{equation}\label{LK}
\|L_K(f)\|_{\ma{H}^k}\leq\|f\|_{\ma{H}^k}.
\end{equation}
\end{prop}

And there are following properties holding for $\ma{L}$, $\ma{L}_{F_*}$ and $\ma{L}_{F_*}^T$:
\begin{lem}\emph{(\cite{Elgindi2019})}We have that
\begin{equation}\label{commu}
L_{K}\left(\mc{L}_{F_*}(f)\right)=\mc{L}\left(L_{K}(f)\right),
\end{equation}and
\begin{equation}
\mc{L}(g)w=gw+z\p_z(gw),
\end{equation}hence
\begin{equation}\label{coer1}
(\mc{L}(g)w,gw)_{L^2}=\f{1}{2}\|gw\|^2_{L^2}.
\end{equation}
\end{lem}

\begin{prop}\emph{(\cite{Elgindi2019})}
Fix $\al<10^{-14}$ and $k\in\mathbb{N}$, there exists $c_k>0$ such that for all $f\in\mc{H}^k$ we have
\begin{equation}\label{coer}
(\mc{L}_{F_*}^T(f),f)_{\mc{H}^k}\geq c_k\|f\|^2_{\mc{H}^k}.
\end{equation}
\end{prop}

\section{Coordinate Transformation, Dynamical Scaling and Linearization}
We work in the framework of \cite{Elgindi2019a}.

\subsection{Vorticity Formation under Cylindrical Coordinate}
Conventionally, it is convenient to consider the blow up behavior of a fuild system in terms of vorticity. Let $\w=\nabla\times u$ and from (\ref{1.1}) we get 
\begin{equation*}
\p_t \omega+(u\cdot\nabla)\w=(\w\cdot\nabla)u+\nabla\times(v\nabla\xi).
\end{equation*}
By Biot-Savart law, if $\td{\psi}$ is the stream field, we have $$-\Delta\td{\psi}=\w,\ \ \ \ \nabla\times\td{\psi}=u.$$
Furthermore, here we consider the axisymmetric case in cylindrical coordinate $(r, \zeta, x^3)$ without swirl, i.e. the velocity field can be written to $u(t,x)=u_r(t, r, x^3)e_r+u_3(t, r, x^3)e_3$,  $v(t,x)=v(t,r,x^3)$, and $u_\zeta\equiv0$, where $r=\sqrt{x_1^2+x_2^2}, \zeta=\arctan (\f{x_2}{x_1})$. Meanwhile, $\w(t,x)=\w_\zeta(t,r,x^3)e_\zeta$.
Moreover, the cylindrical Biot-Savart law reads
\begin{equation*}
\begin{cases}
\p_r(\f{1}{r}\p_r \td{\psi})+\f{1}{r}\p_{33}\td{\psi}=-\w,\\
u_r=\f{1}{r}\p_3\td{\psi},\ \ \ \ u_3=-\f{1}{r}\p_r\td{\psi}.
\end{cases}
\end{equation*}

Now the vorticity form of (\ref{1.1}) turns out to be(we just denote $\w=\w_\zeta$ for simplicity) 
\begin{equation}\label{cyeq}
\begin{cases}
\p_t\w+(u_r\p_r+u_3\p_3)\w=\f{u_r\w}{r}+\p_3v\p_r\xi-\p_rv\p_3\xi,\\
\p_tv+(u_r\p_r+u_3\p_3)v=0,\\
-\p_{rr}\psi-\p_{33}\psi-\f{1}{r}\p_r\psi+\f{\psi}{r^2}=\w,\\
\p_{rr}\xi+\f{\p_r}{r}\xi+\p_{33}\xi=v,
\end{cases}
\end{equation}
where $r\psi=\td{\psi}$.

We further assume that the vorticity satisfies odd condition with respect to $x_3$, i.e. $\w(r,-x_3)=-\w(r,x_3)$, so that $\psi(r,0)=\psi(0,x_3)=0$. Also, we assume charges density satisfies $v(r,-x_3)=v(r,x_3)$ and electrical potential has zero tangential derivative across $x_3$ axis and $x_3=0$ plane, which means $\p_{r}\xi|_{r=0}=0$ and $\p_{x_3}\xi|_{x_3=0}=0$. In addition, we let $\xi\to0$ as $\rho\to\infty$.

\subsection{Spherical Coordinate and Self-similar Transform}

Now we pass to the spherical coordinate. We set $\w(r,x_3)=\W(R,\theta)$, $\psi(r,x_3)=\rho^2\Psi(R,\theta)$, $v(r,x_3)=V(R,\theta)$ and $\xi=\rho^2\Xi(R,\se)$, then we haveW
\begin{equation*}
\p_r\to\f{\cos\theta}{\rho}\al R\p_R-\f{\sin\theta}{\rho}\p_\theta, \ \ \ \ 
\p_3\to\f{\sin\theta}{\rho}\al R\p_R+\f{\cos\theta}{\rho}\p_\theta,
\end{equation*}
and 
\begin{align*}
u_r&=\rho(2\sin\theta\Psi+\al\sin\theta R\p_R\Psi+\cos\theta\p_\theta\Psi),\\
u_3&=\rho(-\f{1}{\cos\se}\Psi-2\cos\theta\Psi-\al\cos\theta R\p_R\Psi+\sin\theta\p_\theta\Psi).
\end{align*}

In this way, (\ref{cyeq}) turns out to be
\begin{equation}\label{newcoor}
\begin{cases}
\p_t \W+\ma{U}(\Psi)\p_\se\W+\ma{V}(\Psi)\al R\p_R\W=\ma{R}(\Psi)\W+\al R\p_R\Xi\p_\se V-\al R\p_RV\p_\se\Xi+2\p_\se V\Xi,\\
\p_t V+\ma{U}(\Psi)\p_\se V+\ma{V}(\Psi)\al R\p_RV=0,\\
\al^2R^2\p_{RR}\Xi+\al(5+\al)R\p_R\Xi+6\Xi+\p_{\se\se}\Xi-\tan\se\p_\se\Xi=V,\\
-\al^2R^2\p_{RR}\Psi-\al(5+\al)R\p_R\Psi-6\Psi+\p_{\se\se}\Psi+\p_\se(\tan\se\Psi)=\W,
\end{cases}
\end{equation}
where $$\ma{U}(\Psi)=-3\Psi-\al R\p_R\Psi,\ \ \ \  \ma{V}(\Psi)=\p_\se\Psi-\tan\se\Psi,$$ and $$\ma{R}(\Psi)=\f{1}{\cos\se}(2\sin\se\Psi+\al\sin\se R\p_R\Psi+\cos\p_\se\Psi).$$

Now we introduce dynamical scaling variables to explore the stability of profiles $(F,\Phi_F)$:
$$z=\f{\mu R}{\la^{1+\delta}},\ \ \ \ \f{ds}{dt}=\f{1}{\la},$$
and 
\begin{equation*}
\W(R,t,\se)=\f{1}{\la}W(z,s,\se),\ \ \ \ \Psi(R,t,\se)=\f{1}{\la}\Phi(z,s,\se),
\end{equation*}
\begin{equation}
\Xi(R,t,\se)=\f{1}{\la}\Pi(z,s,\se),\ \ \ \ V(R,t,\se)=\f{1}{\la}G(z,s,\se).
\end{equation}
So (\ref{newcoor}) becomes
\begin{equation} \label{evo1}
\begin{cases}
\p_s W+\f{\mu_s}{\mu}z\p_zW-\f{\la_s}{\la}\ma{S}_\delta(W)+\ma{U}(\Phi)\p_\se W+\ma{V}(\Phi)\al z\p_zW=\ma{R}(\Phi)W+N_2(\Pi,G),\\
\p_s G+\f{\mu_s}{\mu}z\p_zG-\f{\la_s}{\la}\ma{S}_{\delta}(G)+\ma{U}(\Phi)\p_\se G+\ma{V}(\Phi)\al z\p_zG=0,\\
\al^2z^2\p_{zz}\Pi+\al(\al+5)z\p_z\Pi+6\Pi+\p_{\se\se}\Pi-\tan\se\p_\se\Pi=G,\\
-\al^2z^2\p_{zz}\Phi-\al(\al+5)z\p_z\Phi-6\Phi-\p_{\se\se}\Phi-\tan\se\p_\se\Phi=W,
\end{cases}
\end{equation}where $\ma{S}_\delta(W)=W+(1+\delta)z\p_zW$ and $N_2(\Pi,G)=\al z\p_z\Pi\p_\se G-\al z\p_zG\p_\se\Pi+2\p_\se G\Pi.$

\subsection{Linearization}
In \cite{Elgindi2019} the author shows that there exists a solution to Euler equation of the form $$\W(R,t,\se)=\f{1}{T-t}F\left(\f{R}{(T-t)^{1+\delta}},\se\right),$$where $0<\delta\ll1$ depends on $\al$. 

Recall that $F=F_*+\al^2g$ with $F_*=\f{\al\Gamma}{c}\f{4z}{(1+z)^2}$ and $\|g\|_{\ma{H}^k}\leq C,$ where $C$ independent of $\al$. In particular, 
\begin{equation}\label{LKg}
L_K(g)(0)=0.
\end{equation}

Note that $F$ and $\Phi_F$ satisfy the equations
\begin{equation*}
F+(1+\delta)z\p_zF+\ma{U}(\Phi_F)\p_\se F+\al\ma{V}(\Phi_F)z\p_zF=\ma{R}(\Phi_F)F,
\end{equation*}and 
\begin{equation*}
-\al^2R^2\p_{RR}\Phi_F-\al(\al+5)R\p_R\Phi_F-\p_{\se\se}\Phi_F+\p_\se(\tan\se\Phi_F)-6\Phi_F=F,
\end{equation*}
where $\ma{U}$, $\ma{V}$ and $\ma{R}$ are as above. Now we perform linearization by
$$W=F+\ve,\ \ \Phi=\Phi_F+\Phi_\ve,\ \ G=G+0,$$ and get
\begin{subequations}
\begin{align}
 \label{a}
\p_s\ve+\f{\mu_s}{\mu}z\p_z\ve-(1+\f{\la_s}{\la})\ma{S}_\delta(\ve)+\ma{M}(\ve)=E+N_1(\ve)+N_2(\Pi,G),\\  \label{b}
\p_sG+\f{\mu_s}{\mu}z\p_zG-(1+\f{\la_s}{\la})\ma{S}_\delta(G)+\ma{M}_G(G)=N_3(\ve,G),\\\label{c}
-\al^2z^2\p_{zz}\Phi_\ve-\al(5+\al)z\p_z\Phi_\ve-\p_{\se\se}\Phi_\ve+\p_\se(\tan\se\Phi_\ve)-6\Phi_\ve=\ve,\\ \label{d}
\al^2z^2\p_{zz}\Pi+\al(\al+5)z\p_z\Pi+6\Pi+\p_{\se\se}\Pi-\tan\se\p_\se\Pi=G,
\end{align}
\end{subequations}
where $\ma{M}$ and $\ma{M}_G$ are the linearized operators:
\begin{equation*}
\ma{M}(\ve)=\ma{S}_\delta(\ve)+\ma{U}(\Phi_F)\p_\se\ve+\ma{V}(\Phi_F)\al z\p_z\ve+\ma{U}(\Phi_\ve)\p_\se F+\ma{V}(\Phi_\ve)\al z\p_zF-\ma{R}(\Phi_F)\ve-\ma{R}(\Phi_\ve)F,
\end{equation*}

\begin{equation} \label{MG}
\ma{M}_G(G)=\ma{S}_\delta(G)+\ma{U}(\Phi_F)\p_\se G+\ma{V}(\Phi_F)\al z\p_zG.
\end{equation}

$E$ is the error term:
\begin{equation} \label{err}
E=-\f{\mu_s}{\mu}z\p_zF+(1+\f{\la_s}{\la})\ma{S}_\delta(F).
\end{equation}
And the nonlinear terms are
\begin{equation*}
N_1(\ve)=-\ma{U}(\Phi_\ve)\ve-\al\ma{V}(\Phi_\ve)z\p_z\ve+\ma{R}(\Phi_\ve)\ve,
\end{equation*}
\begin{equation}\label{N2}
N_2(\Pi,G)=\al z\p_z\Pi\p_\se G-\al z\p_zG \p_\se\Pi+2\p_\se G\Pi,
\end{equation}
\begin{equation*}
N_3(\ve,G)=-\ma{U}(\Phi_\ve)\p_\se G-\ma{V}(\Phi_\ve)\al z\p_zG.
\end{equation*}

By elliptic estimate (\ref{el1}), $\ma{U}(\Phi_F)$, $\ma{U}(\Phi_\ve)$ and other similar terms have simple asymtotic forms:

\begin{equation} \label{Ue}
\ma{U}(\Phi_\ve)=-\f{3}{4\al}\sin2\se L_{K}(\ve)+O(1),
\end{equation}
\begin{equation}\label{Ve}
\ma{V}(\Phi_\ve)=\f{1}{4\al}(2\cos2\se-2\sin^2\se)L_{K}(\ve)+O(1),
\end{equation}
\begin{equation*}
\ma{R}(\Phi_\ve)=\f{1}{2\al}L_{K}(\ve)+O(1),
\end{equation*}
where the $O(1)$ terms are bounded by constant.

Recall $F=F_*+\al^2g$ with $$F_*=\f{\Gamma}{c}\f{4\al z}{(1+z)^2},$$
We also have
\begin{equation} \label{UF}
\ma{U}(\Phi_F)=-3\sin2\se\f{1}{1+z}+O(\al),
\end{equation}
\begin{equation} \label{VF}
\ma{V}(\Phi_F)=(2\cos2\se-2\sin^2\se)\f{1}{1+z}+O(\al),
\end{equation}
\begin{equation*}
\ma{R}(\Phi_F)=\f{2}{1+z}+O(\al).
\end{equation*}
Thus, we can rewrite $\ma{M}$ as
\begin{equation} \label{M2}
\ma{M}(\ve)=\ma{L}_{F_*}^T(\ve)+\Gamma(\se)\f{2z^2}{c(1+z)^3}L_K(\f{3}{1+z}\sin2\se\p_\se\ve)(0)+\sqrt{\al}\td{L}(\ve),
\end{equation}
where
$$\td{L}(\ve)=-\f{1}{\sqrt{\al}}\big[\al \ma{V}(F_*)z\p_z\ve+\ma{U}(\Phi_\ve)\p_\se F_*+\al\ma{V}(\Phi_\ve)z\p_zF_*+l.o.t. \big].$$

\section{Coercivity}
In \cite{Elgindi2019a}, the author gives the coericivity of $\ma{M}$:
\begin{prop}\emph{(\cite{Elgindi2019a})}
For all $\ve\in\ma{H}^k$ with $L_{12}(\ve)(0)=0,$ there exists $C_{\ma{M}}>0$ depending only on $k$ such that if $\al<10^{-14}$, we have
\begin{equation*}
(\ma{M}(\ve),\ve)_{\ma{H}^k}\geq C_{\ma{M}}\|\ve\|^2_{\ma{H}^k}.
\end{equation*}
\end{prop}

Then as for $\ma{M}_G$, we have:
\begin{prop}
For all $G\in\ma{H}^k$, there exists $C_{\ma{M}_G}>0$ depending only on $k$ such that if $\al$ is sufficiently small, we have
\begin{equation} \label{coer2}
(\ma{M}_G(G),G)_{\ma{H}^k}\geq C_{\ma{M}_G}\|G\|^2_{\ma{H}^k}.
\end{equation} 
\end{prop}
\proof
By (\ref{MG}), (\ref{UF}), (\ref{VF}), we get
\begin{align*}
\ma{M}_G(G)=&\ma{S}_\delta(G)+\ma{U}(\Phi_F)\p_\se G+\ma{V}(\Phi_F)\al z\p_zG\\
=&G+(1+\delta)z\p_zG-\f{3\sin2\se}{1+z}\p_\se G+\f{\al}{1+z}(2\cos2\se-2\sin^2\se)z\p_z G+O(\al)\\
=&G+z\p_zG-\f{3\sin2\se}{1+z}\p_\se G+l.o.t.\\
=&\ma{L}(G)-\f{3\sin2\se}{1+z}\p_\se G+\f{2G}{1+z}+l.o.t.,
\end{align*}
where $l.o.t.$ stands for terms bounded by $C\al$ in $\ma{H}^k$, which are much smaller than 1 if let $\al\ll1$.

\textbf{Step 1}: $k=0$. 

By the definition of $\ma{M}_G$, (\ref{coer1}) and integration by parts, we have
\begin{align*}
&(\ma{M}_G(G),G\f{w^2}{\sin2\se^\eta})_{L^2}\\
=&\int_{0}^{\infty}\int_{0}^{\pi/2} (\ma{L}G-\f{3}{1+z}\sin2\se\p_\se G+\f{2G}{1+z})G \f{w^2}{\sin2\se^\eta}d\se dz\\
=&\f{1}{2}\|G\f{w}{\sin2\se^{\eta/2}}\|^2_{L^2}+\int_{0}^{\infty}\int_{0}^{\pi/2} \f{3(1-\eta)}{1+z}G^2\cos2\se(\sin2\se)^{-\eta} w^2 d\se dz\\
&+\int_{0}^{\infty}\int_{0}^{\pi/2}\f{2G^2}{1+z}\f{w^2}{\sin2\se^\eta}d\se dz.
\end{align*}
Because $3(1-\eta)\cos2\se=\f{3}{100}\cos2\se<2$, the second term can be absorbed by the third term in the last formula.

So for $k=0$, we have 
\begin{equation}\label{k=0}
(\ma{M}_G(G),G\f{w^2}{\sin2\se^\eta})_{L^2}\geq\f{1}{2}\|G\f{w}{\sin2\se^{\eta/2}}\|^2_{L^2}.
\end{equation}

\textbf{Step 2}: $k=1$.

 We need to estimate the terms involving $D_\se$ and $D_z$.

1. Terms involving $D_\se$:
We have
\begin{align*}
&(D_\se(\ma{M}_G(G)),D_\se G\f{w^2}{\sin2\se^\ga})_{L_2}+O(\al)\\
=&\int_{0}^{\infty}\int_{0}^{\pi/2} (D_\se(\ma{L}G)-\f{3}{1+z}\sin2\se\p_\se D_\se G+\f{2D_\se G}{1+z})D_\se G \f{w^2}{\sin2\se^\ga}d\se dz,\\
=&I_1+I_2+I_3.
\end{align*}

By (\ref{coer1}),
$$I_1=\int_{0}^{\infty}\int_{0}^{\pi/2} (D_\se(\ma{L}G) G \f{w^2}{\sin2\se^\ga}d\se dz=\f{1}{2}\|D_\se G\f{w}{\sin2\se^{\ga/2}}\|^2_{L^2}.$$

For $I_2$ we use integration by parts and get
\begin{align*}
I_2=&\int_{0}^{\infty}\int_{0}^{\pi/2}-\f{3}{1+z}\sin2\se\p_\se D_\se GD_\se G\f{w^2}{sin2\se^\ga}d\se dz\\
=&\int_{0}^{\infty}\int_{0}^{\pi/2} \f{3\al\cos2\se}{20(1+z)}(D_\se G)^2\f{w^2}{\sin2\se^\ga} d\se dz,
\end{align*}
which can be absorbed by $I_3$ as long as $\al$ is small enough.

So we get 
\begin{equation}\label{k=1Dse}
(D_\se(\ma{M}_G(G)),D_\se G\f{w^2}{\sin2\se^\ga})_{L_2}\geq\f{1}{2}\|D_\se G\f{w}{\sin2\se^{\ga/2}}\|^2_{L^2}.
\end{equation}

2. Terms involving $D_z$:
Direct computation gives
\begin{align*}
D_z(\ma{L}G)=&
D_z(G+z\p_zG-\f{2G}{1+z})\\
=&D_zG+z\p_z(D_zG)-\f{2(D_zG)}{1+z}+\f{2zG}{(1+z)^2}\\
=&\ma{L}(D_zG)+\f{2zG}{(1+z)^2},
\end{align*}
and 
\begin{align*}
D_z(\f{2G}{1+z})=\f{2D_zG}{1+z}-\f{2zG}{(1+z)^2}.
\end{align*}

This, together with (\ref{coer1}) and integration by parts, yields
\begin{align*}
&(D_z\ma{M}_G(G),D_zG\f{w^2}{\sin2\se^\eta})_{L^2}+O(\al)\\
=&\int_{0}^{\infty}\int_{0}^{\pi/2}\left[D_z(\ma{L}G)-D_z(\f{3}{1+z}\sin2\se\p_\se G)+D_z\f{2G}{1+z} \right]D_zG   \f{w^2}{\sin2\se^\eta}d\se dz\\
=&\int_{0}^{\infty}\int_{0}^{\pi/2}\left[\ma{L}(D_zG)+\f{3z}{(1+z)^2}\sin2\se\p_\se G-\f{3}{1+z}\sin2\se\p_\se(D_zG)+\f{2D_zG}{1+z}d\se\right]D_zG\f{w^2}{\sin2\se^\eta} dz,\\
=&I_1+I_2+I_3+I_4.
\end{align*}

Still by (\ref{coer1}),  $I_1=\f{1}{2}\|D_zG\f{w}{\sin2\se^{\eta/2}}\|^2_{L^2}$.

And $I_2$ is bounded below by Cauchy-Schwarz
\begin{align*}
I_2=&\int_{0}^{\infty}\int_{0}^{\pi/2}\f{3z}{(1+z)^2}D_\se GD_zG\f{w^2}{\sin2\se^\eta}d\se dz\\
\geq&-\|D_\se G\f{w}{\sin2\se^{\eta/2}}\|_{L^2}\|D_zG\f{w}{\sin2\se^{\eta/2}}\|_{L^2}\\
\geq&-\|D_\se G\f{w}{\sin2\se^{\ga/2}}\|_{L^2}\|D_zG\f{w}{\sin2\se^{\eta/2}}\|_{L^2}\\
\geq&-\f{1}{4}\|D_\se G\f{w}{\sin2\se^{\ga/2}}\|_{L^2}^2-4\|D_zG\f{w}{\sin2\se^{\eta/2}}\|_{L^2}^2,
\end{align*}
where the second inequality holds because $\f{99}{100}=\eta<\ga=1+\f{\al}{10}$ and $\sin2\se^{\eta/2}>\sin2\se^{\ga/2},\se\in[0,\pi/2]$, and the third inequality holds by Young's inequality.

For $I_3$ we use integration by parts again and get
\begin{align*}
I_3=&\int_{0}^{\infty}\int_{0}^{\pi/2}-\f{3}{1+z}\sin2\se\p_\se(D_zG)D_zG\f{w^2}{\sin2\se^\eta}d\se dz\\
=&\int_{0}^{\infty}\int_{0}^{\pi/2}\f{3(1-\eta)\cos2\se}{2(1+z)}(D_zG)^2\f{w^2}{\sin2\se^\eta}d\se dz,
\end{align*}
which can be absorbed by the fourth term for $\f{3}{200}\cos2\se<2$.

So as for terms involving $D_z$ we get 
\begin{equation}\label{k=1Dz}
(D_z\ma{M}_G(G),D_zG\f{w^2}{\sin2\se^\eta})_{L^2}\geq\f{1}{4}\|D_\se G\f{w}{\sin2\se^{\ga/2}}\|_{L^2}^2-4\|D_zG\f{w}{\sin2\se^{\eta/2}}\|^2_{L^2}.
\end{equation}

Then gathering (\ref{k=0}) (\ref{k=1Dse}) and (\ref{k=1Dz}), we have
\begin{align*}
(\ma{M}_G(G),G\f{w^2}{\sin2\se^\eta})_{L_2}+9(D_\se(\ma{M}_G(G)),D_\se G\f{w^2}{\sin2\se^\ga})_{L_2}+(D_z\ma{M}_G(G),D_zG\f{w^2}{\sin2\se^\eta})_{L_2}\\ \geq
\f{1}{2}\|G\f{w}{\sin2\se^{\eta/2}}\|^2_{L^2}+\f{1}{8}\|D_zG\f{w}{\sin2\se^{\eta/2}}\|^2_{L^2}+\f{1}{2}\|D_\se G\f{w}{\sin2\se^{\eta/2}}\|^2_{L^2}.
\end{align*}

\textbf{Step 3}: $k\geq2.$

Recall the definition of $\ma{H}^k$ (\ref{Hk}), we have the following equivalent relation
$$(f,f)_{\ma{H}^k}\approx(f,f)_{\ma{H}^{k-1}}+\|D_\se f\|^2_{\ma{H}^{k-1}}+\|D_zf\|^2_{\ma{H}^{k-1}}\approx \|f\|^2_{\ma{H}^k}.$$

We prove by induction on $k$. Assuming for $k-1$
\begin{align}\label{assu}
(\ma{M}_G(G),G)_{\ma{H}^{k-1}}\geq c_{k-1}\|G\|^2_{\ma{H}^{k-1}},
\end{align}
we prove (\ref{coer2}).

First we note that
\begin{align*}
D_\se(\ma{M}_G(G))=&D_\se(\ma{L}G)-\f{3}{1+z}\sin2\se\p_\se D_\se G+\f{2D_\se G}{1+z}\\
=&\ma{M}_G(D_\se G),
\end{align*}
so by (\ref{assu}) we have 
\begin{equation}\label{Dsk}
(D_\se(\ma{M}_G(G)),D_\se(G))_{\ma{H}^{k-1}}\geq c_{k-1}\|D_\se G\|^2_{\ma{H}^{k-1}}.
\end{equation}

And note that
\begin{align*}
D_z(\ma{M}_G(G))=&D_z(\ma{L}G)-D_z(\f{3}{1+z}\sin2\se\p_\se G)+D_z\f{2G}{1+z}\\
=&\ma{L}(D_zG)+\f{3z}{(1+z)^2}\sin2\se\p_\se G-\f{3}{1+z}\sin2\se\p_\se(D_zG)+\f{2D_zG}{1+z}\\
=&\ma{M}_G(D_zG)+\f{3z}{(1+z)^2}\sin2\se\p_\se G,
\end{align*}
by (\ref{assu}) and Young's inequality we get
\begin{align*}
&(D_z(\ma{M}_G(G)),D_z(G))_{\ma{H}^{k-1}}\\
=&(\ma{M}_G(D_zG),D_zG)_{\ma{H}^{k-1}}+(\f{3z}{(1+z)^2}\sin2\se\p_\se G,D_zG)_{\ma{H}^{k-1}}\\
\geq&c_{k-1}\|D_zG\|^2_{\ma{H}^{k-1}}-\f{c_{k-1}}{2}\|D_zG\|^2_{\ma{H}^{k-1}}-C_{k-1}\|D_\se G\|^2_{\ma{H}^{k-1}},
\end{align*}
that is 
\begin{equation}\label{Dzk}
(D_z(\ma{M}_G(G)),D_z(G))_{\ma{H}^{k-1}}\geq\f{c_{k-1}}{2}\|D_zG\|^2_{\ma{H}^{k-1}}-C_{k-1}\|D_\se G\|^2_{\ma{H}^{k-1}}.
\end{equation}

Gathering the (\ref{assu}), (\ref{Dsk}) and (\ref{Dzk}), there exist $c_k>0$ such that 
\begin{equation*}
(\ma{M}_G(G),G)_{\ma{H}^{k}}\geq c_k\|G\|^2_{\ma{H}^k}.
\end{equation*}
\qed

\section{Analysis of Variant Spherical Laplacian}
Now we consider the elliptic equations with boundary conditions assumed in Section 2. First we consider
\begin{equation}\label{ellp1}
\begin{cases}
-\al^2z^2\p_{zz}\Psi-\al(5+
\al) z\p_z\Psi-\p_{\se\se}\Psi+\p_\se(\tan\se\Psi)-6\Psi=F,\\
\Psi(z,0)=0,\ \ \Psi(z,\f{\pi}{2})=0,\ \ \Psi(\infty,\se)=0,
\end{cases}
\end{equation}which arises from Biot-Savart law. In \cite{Elgindi2019} the author gives the following estimate on (\ref{ellp1}):
\begin{prop}
Let $\al>0$ and assume $F\in\ma{H}^k$. Let $\Psi$ be the unique solution to (\ref{ellp1}) which vanishes at $\se=0,\f{\pi}{2}$ as $z\to \infty$. Then
\begin{equation}\label{el1}
\al^2\|D^2_z\Psi\|_{\ma{H}^k}+\al\|D_z\Psi\|_{\ma{H}^k}+\|\p_{\se\se}(\Psi-\f{1}{4\al}\sin2\se L_{K}(F))\|_{\ma{H}^k}\leq C_k\|F\|_{\ma{H}^k},
\end{equation}
where $L_{K}(F)=\int_{z}^{\infty}\int_{0}^{\pi/2}\f{F(\rho,\se)K(\se)}{z}d\rho d\se$.
\end{prop}

Next, we consider
\begin{equation}\label{ellp2}
\begin{cases}
\al^2z^2\p_{zz}\Pi+\al(5+\al)z\p_z\Pi+\p_{\se\se}\Pi-\tan\se\p_\se\Pi+6\Pi=G,\\
\p_\se\Pi(z,0)=0,\ \ \p_\se\Pi(z,\f{\pi}{2})=0,\ \ \Pi(\infty,\se)=0,
\end{cases}
\end{equation}which is derived from variant spherical Laplacian  corresponding to electrical potential. We want to find a special solution of (\ref{ellp2}).

First we consider the part relevant to $z$:
\begin{equation}
\begin{cases}
\al^2z^2\p_{zz}\Pi_1+\al(5+\al)z\p_z\Pi_1+\Pi_1=0,\\
\Pi_1(\infty)=0,
\end{cases}
\end{equation}
which is an ODE of Euler type, and $\Pi_1(z)=c_1z^{\f{\sqrt{21}-5}{2\al}}$ is a special solution, where $c_1$ is under determined.
Secondly, by Sturm-Liouville theory, direct computation shows that $\Pi_2(\se)=\cos^2\se-\f{2}{3}$ is a solution of
\begin{equation*}
\begin{cases}
\p_{\se\se}\Pi_2-\tan\se\p_\se\Pi_2+6\Pi_2=0,\\
\p_\se\Pi_2(0)=\p_\se\Pi(\f{\pi}{2})=0.
\end{cases}
\end{equation*}

Then, we construct $\Pi_0$ as
\begin{equation*}
\Pi_0=\td{\Pi}(s)\Pi_1(z)\Pi_2(\se),
\end{equation*}
where $\td{\Pi}(s)$ is under determined.
Obviously, $\Pi_0$ satisfies
$$\al^2z^2\p_{zz}\Pi_0+\al(5+\al)z\p_z\Pi_0+\p_{\se\se}\Pi_0-\tan\se\p_\se\Pi_0+6\Pi_0=-\Pi_0.$$

However, $\Pi_0\notin\ma{H}^k$ because $z^{\f{\sqrt{21}-5}{2\al}}$ grows too fast near origin, so we let 
\begin{equation} \label{Pi}
\Pi=\chi\Pi_0,
\end{equation}where
\begin{equation*}
\chi(z)=
\begin{cases}0\ \ \ z<\f{1}{2},\\
1\ \ \ z>1.
\end{cases}
\end{equation*}

We also have the following properties hold for $\Pi$,
\begin{prop}
Let $P=\al^2z^2\p_{zz}+\al(5+\al)z\p_z+\p_{\se\se}-\tan\se\p_\se+6Id$, and $G=P\Pi$, then the following hold:
\begin{itemize}
\item$\Pi\in\ma{H}^k$. In particular, $\|\Pi\|_{\ma{H}^k}=C\td{\Pi}(s)$ for some $C$;
\item$|G|$=$|P\Pi|\leq C|\Pi|$, where $C$ is independent of $z$, so $(G,G)_{\ma{H}^{k}}=\td{\Pi}^2(s)$ for suitable $c_1$;
\item$\p_\se\Pi(0)=\p_\se\Pi(\f{\pi}{2})=\p_zG(0,\se)=0$.
\end{itemize}
\end{prop}
\proof
We have
$$\|\Pi\|_{\ma{H}^k}=\td{\Pi}(s)\|\chi(z)z^{\f{\sqrt{21}-5}{2\al}}(\cos^2\se-\f{2}{3})\|_{\ma{H}^k}=C\td{\Pi}(s),$$ so the first assertion holds.
Meanwhile,
\begin{align*}
P\Pi&=P(\chi\Pi_0)\\
&=(\al^2z^2\p_{zz}+\al(5+\al)z\p_z)(\chi\Pi_0)\\
&=\al^2z^2(\p_{zz}\chi\Pi_0+2\p_z\chi\p_z\Pi_0)+\al(5+\al)z\p_z\chi\Pi_0+\chi\Pi_0.
\end{align*}
It is easy to check $|\al^2z^2\p_{zz}\chi|+|\al(5+
\al)z\p_z\chi| \leq C|\chi|$ and $|\al^2z^2\p_z\chi\p_z\Pi_0|\leq C|\chi\Pi_0|$, so the second assertion holds.
The third assertion is obvious.
\section{Modulation}
In this section, for any $s>0$, we choose proper $\la$ and $\mu$ such that  
\begin{equation}\label{LK0}L_{K}(\ve)(0)=0,\end{equation}
to impose (\ref{LK}) to $\ve(z,s,\se)$. And note that if $\p_z\ve(0,s,\se)=0$ then 

\begin{equation}\label{LKz}
L_{K}(z\p_z\ve)(0)=0,\ \ \ \  s>0.
\end{equation}

\begin{prop}
For any $\se\in[0,\pi/2]$, $s>0$, in order to keep $L_{K}(\ve)(0)=0$ and $\p_z\ve(0,s,\se)=0$, it suffices to impose $\la$ and $\mu$ to satisfy
\begin{equation}\label{modu1}
\al(\f{\la_s}{\la}+1)-3L_{K}(\f{\sin2\se}{1+z}\p_\se\ve)(0)+\sqrt{\al}L_K(\td{L}(\ve))(0)+L_K(N_1(\ve))(0)+L_K(N_2(\Pi,G))(0)=0,
\end{equation}
\begin{equation*}
\la(0)=1,
\end{equation*}
 and 
\begin{equation}\label{modu2}
\f{\mu_s}{\mu}=(2+\delta)(\f{\la_s}{\la}+1).
\end{equation} 
In particular, we have the following estimates:
\begin{equation}\label{mo1}
|\al(\f{\la_s}{\la}+1)-3L_{K}(\f{\sin2\se}{1+z}\p_\se\ve)(0)|\les\sqrt{\al}\|\ve\|_{\ma{H}^k}+\al\|\ve\|^2_{\ma{H}^k}+\td{\Pi}^2(s),
\end{equation}

\end{prop}

and
\begin{align}\label{roughmodu}
|\f{\la_s}{\la}+1|\les\f{1}{\al}\|\ve\|_{\ma{H}^k}.
\end{align}

\proof
At first, since $\p_zN_2(\Pi,G)(0)=0$, the proof of (\ref{modu2}) is just a repeatation of that in \cite{Elgindi2019a}, and $\p_z\ve(0,s,\se)=0$ follows.

It remains to keep $L_K(\ve)(0)=0$. We impose $L_{K}$ on the first equation of (\ref{a}) and let $z=0$ to get
\begin{align*}
\p_sL_{K}(\ve)(0)+\f{\mu_s}{\mu}L_{K}(z\p_z\ve)(0)-(\f{\la_s}{\la}+1)L_K(\ma{S}_\delta(\ve))(0)=L_K(E)(0)-L_K(\ma{M}\ve)(0)\\+L_K(N_1(\ve))(0)+L_K(N_2(\Pi,G))(0).
\end{align*}
By (\ref{zpz}), (\ref{err}) and (\ref{LKz}), as well as $F=F_*+\al^2g$, we have
\begin{align*}
L_K(E)(0)=&L_K(-\f{\mu_s}{\mu}z\p_zF+(1+\f{\la_s}{\la})\ma{S}_\delta(F))(0)\\
=&L_K((1+\f{\la_s}{\la})F)(0)\\
=&L_K((1+\f{\la_s}{\la})(F_*+\al^2g))(0)\\
=&4\al((1+\f{\la_s}{\la})).
\end{align*}

Meanwhile, (\ref{M2}) and (\ref{commu}) yield
\begin{align*}
L_K(\ma{M}\ve)(0)=&L_K(\ma{L}_{F_*}^T(\ve)+\Gamma(\se)\f{2z^2}{c(1+z)^3}L_K(\f{3}{1+z}\sin2\se\p_\se\ve)(0)+\sqrt{\al}\td{L}(\ve))(0)\\
=&L_K(\ma{L}_{F_*}\ve)(0)-3L_K(\f{\sin2\se}{1+z}\p_\se\ve)(0)+\sqrt{\al}L_K(\td{L}\ve)(0)\\
=&\ma{L}(L_K(\ve))(0)-3L_K(\f{\sin2\se}{1+z}\p_\se\ve)(0)+\sqrt{\al}L_K(\td{L}\ve)(0),
\end{align*}

Then,
\begin{align*}
\p_sL_{K}(\ve)(0)-(\f{\la_s}{\la}+1)L_K(\ve)(0)=4\al((1+\f{\la_s}{\la}))+\ma{L}(L_K(\ve))(0)-3L_K(\f{\sin2\se}{1+z}\p_\se\ve)(0)\\+\sqrt{\al}L_K(\td{L}\ve)(0)+L_K(N_1(\ve))(0)+L_K(N_2(\Pi,G))(0).
\end{align*}
So (\ref{modu1}) implies 
\begin{align*}
\p_sL_{K}(\ve)(0)=&(\f{\la_s}{\la}+1)L_K(\ve)(0)+\ma{L}(L_K(\ve))(0)\\
=&(\f{\la_s}{\la}+1)L_K(\ve)(0)+L_K(\ve)(0)+(z\p_zL_K(\ve))(0)-(\f{2}{1+z}L_K(\ve))(0)\\
=&\f{\la_s}{\la}L_K(\ve)(0)-\int_{0}^{\pi/2}\ve(0,s,\se) K(\se)d\se\\
=&\f{\la_s}{\la}L_K(\ve)(0),
\end{align*}
where we used (\ref{zpz}) in the third equality and $\ve(0,s,\se)=0$ in the last equality. Then we have
$$L_K(\ve)(0,s)=\f{L_K(\ve)(0,0)}{\la(0)}\la(s),$$
so $L_K(\ve)(0)$ will vanish identically once $L_K(\ve)(0)|_{s=0}=0$ and $\la(0)=1$.

By (\ref{modu1}) and \cite{Elgindi2019a}, it easy to get 
$$|4\al(\f{\la_s}{\la}+1)-3L_{K}(\f{\sin2\se}{1+z}\p_\se\ve)(0)|\les\sqrt{\al}\|\ve\|_{\ma{H}^k}+\al\|\ve\|^2_{\ma{H}^k}+\td{\Pi}^2(s)$$
and
$$|\f{\la_s}{\la}+1|\les\f{1}{\al}|L_K(\f{\sin2\se}{1+z}\p_\se\ve)(0)|\les\f{1}{\al}\|\ve\|_{\ma{H}^k}.$$

\qed

\section{Energy Estimate}
In this section we derive energy esitmate for (\ref{a}) and (\ref{b}).

First we define the energy
\begin{align*}
\ma{E}(s)=(\ve,\ve)_{\ma{H}^k}+\td{\Pi}^2(s).
\end{align*}

\begin{prop}If $\al\ll1$, we have:
\begin{equation}\label{energy}
\f{d}{ds}\ma{E}\leq-c\ma{E}+\f{C}{\al^{3/2}}\ma{E}^{3/2}.
\end{equation}for some $c$, $C>0$.
\end{prop}
\proof

First, we check the boundness of
$\f{d}{ds}(\ve,\ve)_{\ma{H}^k}$. By (\ref{a}), we have
\begin{align*}
\f{1}{2}\f{d}{ds}(\ve,\ve)_{\ma{H}^k}\leq-(\ma{M}\ve,\ve)_{\ma{H}^k}+(E,\ve)_{\ma{H}^k}+|\f{\mu_s}{\mu}|(z\p_z\ve,\ve)_{\ma{H}^k}+|(1+\f{\la_s}{\la})|(\ma{S}_\delta\ve,\ve)_{\ma{H}^k}\\+(N_1(\ve),\ve)_{\ma{H}^k}+(N_2(\Pi,G),\ve)_{\ma{H}^k}.
\end{align*}

As in \cite{Elgindi2019a}, we have
$$-(\ma{M}\ve,\ve)_{\ma{H}^k}+(E,\ve)_{\ma{H}^k}\leq -c\|\ve\|^2_{\ma{H}^k}+\|\ve\|_{\ma{H}^k}(\sqrt{\al}\|\ve\|_{\ma{H}^k}+\al\|\ve\|_{\ma{H}^k}^2)+\al\|g\|_{\ma{H}^k}\|\ve\|^2_{\ma{H}^k},$$
$$|\f{\mu_s}{\mu}|(z\p_z\ve,\ve)_{\ma{H}^k}\leq\f{C}{\al^{3/2}}\|\ve\|^3_{\ma{H}^k},$$
$$|\f{\la_s}{\la}+1||(\ma{S}_\delta(\ve),\ve)_{\ma{H}^k}|\leq\f{C}{\al^{3/2}}\|\ve\|^3_{\ma{H}^k},$$and
$$(N_1(\ve),\ve)_{\ma{H}^k}\leq\f{C}{\al^{3/2}}\|\ve\|^3_{\ma{H}^k}.$$

As for $(N_2(\Pi,G),\ve)_{\ma{H}^k}$, since $\Pi=\td{\Pi}(s)\chi(z)z^{\f{\sqrt{21}-5}{2\al}}(\cos^2\se-\f{2}{3})$, $|P\Pi|\leq C|\Pi|$, and $\|\Pi\|_{\ma{H}^{k}}=C\td{\Pi}(s)$, we get 
$$\|\p_\se G\|_{\ma{H}^{k}}\leq C\|\p_\se\Pi\|_{\ma{H}^{k}}\leq C\|\Pi\|_{\ma{H}^{k}}\leq C\td{\Pi}(s).$$

Then, by the product lemmas in the Appendix of \cite{Elgindi2019a}, we have
$$(\al z\p_z\Pi\p_\se G,\ve)_{\ma{H}^k}\leq C\al^{1/2}\|\Pi\|_{\ma{H}^{k}}\|G\|_{\ma{H}^{k}}\|\ve\|_{\ma{H}^k}\leq C\al^{1/2}\td{\Pi}^2(s)\|\ve\|_{\ma{H}^k},$$
$$(\al z\p_z G\p_\se \Pi,\ve)_{\ma{H}^k}\leq C\al^{1/2}\|\Pi\|_{\ma{H}^{k}}\|G\|_{\ma{H}^{k}}\|\ve\|_{\ma{H}^k}\leq C\al^{1/2}\td{\Pi}^2(s)\|\ve\|_{\ma{H}^k},$$
$$(\p_\se G\Pi,\ve)_{\ma{H}^k}\leq\f{C}{\al^{1/2}}\|\Pi\|_{\ma{H}^{k}}\|G\|_{\ma{H}^{k}}\|\ve\|_{\ma{H}^k}\leq\f{C}{\al^{1/2}}\td{\Pi}^2(s)\|\ve\|_{\ma{H}^k}.$$
So by (\ref{N2}), $$(N_2(\Pi,G),\ve)_{\ma{H}^k}\leq(\f{C}{\al^{1/2}}+C\al^{1/2})\td{\Pi}^2(s)\|\ve\|_{\ma{H}^k}.$$

Now we check the boundness of $\f{d}{ds}\td{\Pi}^2(s)=\f{d}{ds}(G,G)_{\ma{H}^{k}}.$ By (\ref{b}), we have
$$\f{1}{2}\f{d}{ds}(G,G)_{\ma{H}^{k}} \leq-(\ma{M}_G(G),G)_{\ma{H}^{k}}+|\f{\mu_s}{\mu}|(z\p_z G,G)_{\ma{H}^{k}}+|(1+\f{\la_s}{\la})|(\ma{S}_\delta G,G)_{\ma{H}^{k}}+(N_3(\ve,G),G)_{\ma{H}^{k}}.$$

By (\ref{coer2}), we have
$$-(\ma{M}_G(G),G)_{\ma{H}^{k}}\leq-c\|G\|^2_{\ma{H}^{k}}=-c\td{\Pi}^2(s),$$
and meanwhile,
$$|\f{\mu_s}{\mu}|(z\p_z G,G)_{\ma{H}^{k}}\leq\f{C}{\al^{3/2}}\|G\|^2_{\ma{H}^{k}}\|\ve\|_{\ma{H}^{k}},$$
$$|(1+\f{\la_s}{\la})|(\ma{S}_\delta G,G)_{\ma{H}^{k}}\leq\f{C}{\al^{3/2}}\|G\|^2_{\ma{H}^{k}}\|\ve\|_{\ma{H}^{k}}.$$

And for the nonlinear term $(N_3(\ve),G)_{\ma{H}^{k}}$, because 
$$\|z\p_z G\|_{\ma{H}^{k}}\leq C\|z\p_z\Pi\|_{\ma{H}^{k}}\leq C\|\Pi\|_{\ma{H}^{k}}\leq C\td{\Pi}(s),$$ as well as (\ref{Ue}), (\ref{LK}), (\ref{LK0}) and product lemma in \cite{Elgindi2019a}, we have
$$(\ma{U}(\Phi_\ve)\p_\se G,G)_{\ma{H}^{k}}\leq\f{C}{\al}\|G\|^2_{\ma{H}^{k}}\|L_K(\ve)\|_{\ma{H}^k}\leq\f{C}{\al}\td{\Pi}^2(s)\|\ve\|_{\ma{H}^k},$$
and
$$(\ma{V}(\Phi_\ve)\al z\p_zG,G)_{\ma{H}^{k+1}}\leq C\td{\Pi}^2(s)\|\ve\|_{\ma{H}^k}.$$

Gathering preceding estimates together, (\ref{energy}) follows.
\qed

Now we get another version of Theorem 1.1:
\begin{thm}
For $k\geq4$, there exists $0<\al_0\ll1$ and $\nu,\kappa_0>0$, such that if $\Pi$ is defined as (\ref{Pi})
, and for every initial $\ma{E}(\ve_0,\td{\Pi}(s_0))$ $<\nu_0\al^{3/2}$ and $L_K(\ve_0)=0$, there exists a unique solution to (\ref{a}) so that:
$$
|\f{\mu_s}{\mu}|+|\f{\la_s}{\la}+1|+\ma{E}(\ve,\td{\Pi})(s)\leq C\ma{E}(\ve_0,\td{\Pi}(s_0))e^{-\kappa s}.
$$
\end{thm}
\proof
By (\ref{energy}), (\ref{modu2}) and (\ref{roughmodu}), modulation  and Gronwall's inequality, the rapid dacay of $|\f{\mu_s}{\mu}|$ and $|\f{\la_s}{\la}+1|$, as well as energy estimate. And the existence and uniqueness can be obtained by weak compactness and the decay of energy similarly in  \cite{Elgindi2019}.


\begin{thebibliography}{10}
	
	\bibitem{Balbuena2004}
	Perla~B Balbuena and Yixuan Wang.
	\newblock {\em Lithium-Ion Batteries}.
	\newblock Imperial College Press and Distributed by World Scientific Publishing
	CO., 2004.
	
	\bibitem{Bazant2004}
	Martin~Z. Bazant and Todd~M. Squires.
	\newblock Induced-charge electrokinetic phenomena: Theory and microfluidic
	applications.
	\newblock {\em Physical Review Letters}, 92(6), 2004.
	
	\bibitem{Biler2000}
	Piotr Biler and Jean Dolbeault.
	\newblock Long time behavior of solutions to nernst-planck and debye-hückel
	drift-diffusion systems.
	\newblock {\em Annales Henri Poincar{\'{e}}}, 1(3):461--472, 2000.
	
	\bibitem{Bothe2014}
	Dieter Bothe, Andr{\'{e}} Fischer, and Jürgen Saal.
	\newblock Global well-posedness and stability of electrokinetic flows.
	\newblock {\em {SIAM} Journal on Mathematical Analysis}, 46(2):1263--1316,
	2014.
	
	\bibitem{Buckmaster2019}
	Tristan Buckmaster, Steve Shkoller, and Vlad Vicol.
	\newblock Formation of point shocks for 3d compressible euler.
	\newblock {\em arXiv:1912.04429}.
	
	\bibitem{Chu2005}
	Kevin~T. Chu and Martin~Z. Bazant.
	\newblock Electrochemical thin films at and above the classical limiting
	current.
	\newblock {\em {SIAM} Journal on Applied Mathematics}, 65(5):1485--1505, 2005.
	
	\bibitem{Collot2018a}
	Charles Collot, Tej-Eddine Ghoul, Slim Ibrahim, and Nader Masmoudi.
	\newblock On singularity formation for the two dimensional unsteady prandtl's
	system.
	\newblock {\em arXiv:1808.05967}.
	
	\bibitem{Constantin2018}
	Peter Constantin and Mihaela Ignatova.
	\newblock On the
	nernst{\textendash}planck{\textendash}navier{\textendash}stokes system.
	\newblock {\em Archive for Rational Mechanics and Analysis}, 232(3):1379--1428,
	2018.
	
	\bibitem{Deng2011}
	Chao Deng, Jihong Zhao, and Shangbin Cui.
	\newblock Well-posedness for the
	navier{\textendash}stokes{\textendash}nernst{\textendash}planck{\textendash}poisson
	system in triebel{\textendash}lizorkin space and besov space with negative
	indices.
	\newblock {\em Journal of Mathematical Analysis and Applications},
	377(1):392--405, 2011.
	
	\bibitem{Elgindi2019}
	Tarek~M. Elgindi.
	\newblock Finite-time singularity formation for $c^{1,\al}$ solutions to the
	incompressible euler equations on $\mathbb{R}^3$.
	\newblock {\em arXiv:1904.04795}.
	
	\bibitem{Elgindi2019a}
	Tarek~M. Elgindi, Tej-Eddine Ghoul, and Nader Masmoudi.
	\newblock On the stability of self-similar blow-up for $c^{1,\al}$ solutions to
	the incompressible euler equations on $\mathbb{R}^3$.
	\newblock {\em arXiv:1910.14071}.
	
	\bibitem{Enikov2000}
	Eniko~T. Enikov and Bradley~J. Nelson.
	\newblock Electrotransport and deformation model of ion exchange membrane-based
	actuators.
	\newblock In Yoseph Bar-Cohen, editor, {\em Smart Structures and Materials
		2000: Electroactive Polymer Actuators and Devices ({EAPAD})}. {SPIE}, 2000.
	
	\bibitem{Enikov2005}
	Eniko~T. Enikov and Geon~S. Seo.
	\newblock Analysis of water and proton fluxes in ion-exchange
	polymer{\textendash}metal composite ({IPMC}) actuators subjected to large
	external potentials.
	\newblock {\em Sensors and Actuators A: Physical}, 122(2):264--272, a 2005.
	
	\bibitem{Fischer2016}
	Andr{\'{e}} Fischer and Jürgen Saal.
	\newblock Global weak solutions in three space dimensions for electrokinetic
	flow processes.
	\newblock {\em Journal of Evolution Equations}, 17(1):309--333, 2016.
	
	\bibitem{Jerome2002}
	Joseph~W. Jerome.
	\newblock Analytical approaches to charges transport in a moving medium.
	\newblock {\em Transport Theory and Statistical Physics}, 31(4-6):333--366,
	2002.
	
	\bibitem{Jerome2008}
	Joseph~W. Jerome, Bice Chini, Massimo Longaretti, and Riccardo Sacco.
	\newblock Computational modeling and simulation of complex systems in
	bio-electronics.
	\newblock {\em Journal of Computational Electronics}, 7(1):10--13, 2008.
	
	\bibitem{Jerome2009}
	Joseph~W. Jerome and Riccardo Sacco.
	\newblock Global weak solutions for an incompressible charged fluid with
	multi-scale couplings: Initial{\textendash}boundary-value problem.
	\newblock {\em Nonlinear Analysis: Theory, Methods {\&} Applications},
	71(12):e2487--e2497, 2009.
	
	\bibitem{Liu2019}
	Qiao Liu.
	\newblock The 3d nonlinear dissipative system modeling electro-diffusion with
	blow-up in one direction.
	\newblock {\em Communications in Mathematical Sciences}, 17(1):131--147, 2019.
	
	\bibitem{Longaretti2008}
	Massimo Longaretti, Bice Chini, Joseph~W. Jerome, and Riccardo Sacco.
	\newblock Electrochemical modeling and characterization of voltage operated
	channels in nano-bio-electronics.
	\newblock {\em Sensor Letters}, 6(1):49--56, 2008.
	
	\bibitem{Merle2019}
	Frank Merle, Pierre Raphael, Igor Rodnianski, and Jeremie Szeftel.
	\newblock On smooth self similar solutions to the compressible euler equations.
	\newblock {\em arXiv:1912.10998}.
	
	\bibitem{Probstein2003}
	Ronald~F. Probstein.
	\newblock {\em Physicochemical Hydrodynamics 2e P}.
	\newblock John Wiley \& Sons, 2003.
	
	\bibitem{Rubinstein1990}
	Isaak Rubinstein.
	\newblock {\em Electro-diffusion of ions}, volume~11 of {\em SIAM Studies in
		Applied Mathematics}.
	\newblock Society for Industrial and Applied Mathematics (SIAM), Philadelphia,
	PA, 1990.
	
	\bibitem{Ryham2009}
	Rolf~J. Ryham.
	\newblock Existence, uniqueness, regularity and long-term behavior for
	dissipative systems modeling electrohydrodynamics.
	\newblock {\em arXiv:0910.4973}.
	
	\bibitem{SCHMUCK2009}
	Markus Schmuck.
	\newblock Analysis of the navier-stokes-nernst-planck-poisson system.
	\newblock {\em Mathematical Models and Methods in Applied Sciences},
	19(06):993--1014, 2009.
	
	\bibitem{Shahinpoor2004}
	Mohsen Shahinpoor and Kwang~J Kim.
	\newblock Ionic polymer{\textendash}metal composites: {III}. modeling and
	simulation as biomimetic sensors, actuators, transducers, and artificial
	muscles.
	\newblock {\em Smart Materials and Structures}, 13(6):1362--1388, 2004.
	
	\bibitem{Zhao2016}
	Jihong Zhao and Meng Bai.
	\newblock Blow-up criteria for the three dimensional nonlinear dissipative
	system modeling electro-hydrodynamics.
	\newblock {\em Nonlinear Analysis: Real World Applications}, 31:210--226, 2016.
	
	\bibitem{Zhao2010}
	Jihong Zhao, Chao Deng, and Shangbin Cui.
	\newblock Global well-posedness of a dissipative system arising in
	electrohydrodynamics in negative-order besov spaces.
	\newblock {\em Journal of Mathematical Physics}, 51(9):093101, 2010.
	
\end{thebibliography}
\end{document}